\documentclass[11pt]{article}
\usepackage[utf8]{inputenc}
\usepackage[
papersize={5.5in, 8.5in},
left=0.35in,
right=0.35in,
top=0.35in,
bottom=0.5in]{geometry}
\usepackage{moresize}
\usepackage{array}
\usepackage{graphicx}
\usepackage{amsthm}
\usepackage{rotating}
\usepackage{amsfonts,amssymb,amsmath}
\usepackage[usenames]{color}
\usepackage[T1]{fontenc}

\usepackage{array}
\usepackage{float}
\usepackage{booktabs}
\usepackage[dvipsnames,usenames,table,xcdraw]{xcolor}

\newcolumntype{P}[1]{>{\centering\arraybackslash}p{#1}}

\newcommand{\m}{\bf}
\newcommand{\comment}[1]{}

\graphicspath{ {images/} }

\title{\bf Tiling the plane with hexagons: improved separations for $k$-colourings}
\author{\normalsize {Aubrey D.N.J. de Grey; \textcolor[rgb]{0.5,0.5,0.5}{Jaan Parts (jaan\_parts@mail.ru)} }}
\date{\normalsize {Mountain View, California, USA; \textcolor[rgb]{0.5,0.5,0.5}{Kazan, Tatarstan, Russia}}}

\begin{document}

\maketitle

\pagestyle{empty}
\thispagestyle{empty}

\begin{abstract}
It has been common knowledge since 1950 that seven colours can be assigned to tiles of an infinite 
honeycomb 
with cells of unit diameter
such that no two tiles of the same colour are closer than $d(7)=\frac{\sqrt{7}}{2}$ apart. Various authors have described tilings using $k>7$ colours, giving corresponding values for $d(k)$, but it is generally unknown whether these are the largest possible for a given $k$. Here, for many $k$, we describe tilings with larger values of $d(k)$ than previously reported. 
\end{abstract}

\section{Background}

\begin{figure}[!b]
\centering
\includegraphics[scale=0.4]{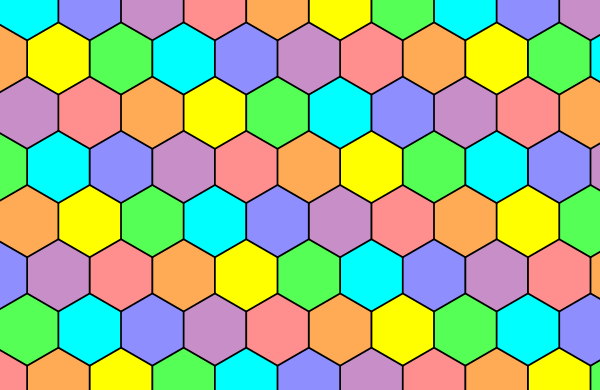} 
\caption{The classical $7$-tiling of the plane with regular hexagons.}
\label{f1}
\end{figure}

In 1950, Isbell observed \cite{soi} that a tiling of the plane using regular hexagons of unit diameter can be $7$-coloured such that tiles of the same colour are at least $d(7)=\frac{\sqrt{7}}{2}$ apart (see Fig.~\ref{f1}). If we allow more colours $k$, this minimum separation $d(k)$ increases. In extension of earlier work, Chybowska-Sokół et al. recently reported \cite{chy} values of $d(k)$ for various $k>7$. Their analysis, like all prior work of which we are aware 
(with a few rare exceptions),
was basically restricted to \textit{regular hexagons}. With such a tiling, if and only if $k$ is a so-called \textit{L\"{o}schian number} (that is, $k=a^2+ab+b^2$ for some integers $a>0, a\ge b\ge 0$), colours can be assigned such that the tiles of any given colour lie on a (regular) 
hexagonal sublattice.
L\"{o}schian numbers are quite common: the first few are $k\in\{1,3,4,7,9,12,13,16,19,21,\dots\}$. In such cases, 
$d(k)$ is $\frac{\sqrt{3}}{2}(a-1)$ when $b=0$, otherwise $\sqrt{(3a+3b-4)^2+3(a-b)^2}/4$. 

It is left to the reader to verify that, in any tiling with regular hexagons, the smallest $k' > k$ that can have all same-coloured pairs of tiles more than $d(k)$ apart is always also a L\"{o}schian number. When multiple choices of $a,b$ give the same $k\in \{49, 91, 133, 147, 169,\dots\}$, the largest $d(k)$ arises from the largest value of $b$.

However, regular hexagons are not the only ones that can tile the plane, even if we restrict ourselves to the case where all tiles are \textit{identical} (or, more correctly, \textit{congruent}). We wondered whether tilings with non-regular hexagons might achieve improved bounds, i.e. larger values for some $d(k)$. Here we report that they indeed do. In fact, we now have 
$d(k) > d(k-3)$ for all $k\le 175$. 

\section{Preliminaries}

\paragraph{Problem statement and main restrictions.}

Without loss of generality, we restrict ourselves to tiles of a \textit{unit diameter}\footnote{
Sometimes it is more correct to use the term \textit{unit width}, referring to the existence of figures with a constant width, such as the Reuleaux triangle and pentagon, which do not fit into a circle of unit diameter. In fact, we require that the distance between any pair of points on a tile be less than one.
}.
Our task is to obtain the maximum possible distance $d(k)$ between tiles of the same colour for each $k$.

We restrict the area of study to the case of periodic tiling of so-called \textit{lattice-sublattice} scheme, 
where
i) the plane is partitioned into $k$ congruent sets of tiles, obtained one from the other by translation and reflection (so for analysis it is enough to consider tiles of one chosen colour),
ii) all tiles have the same shape.

Until now, we have implicitly assumed that $d(k)$ is a non-decreasing function of the number of colours $k$. The restrictions i) and ii) introduced have a slightly unexpected effect: 
some $k$ can give a local decrease in $d(k)$.

\paragraph{A “hierarchy of irregularity” of hexagons that tile the plane.}
From the point of view of maximising $d(k)$, convex hexagons and pentagons are of main interest. 
A full taxonomy of tilings of the plane with identical (congruent) convex hexagons has been described by Gardner \cite{gar}. Also, Rao recently showed \cite{rao} that the known classification of tilings using identical pentagons 
is complete.

Here we restrict ourselves to a narrow subset of hexagons, namely those with all three diagonals being the same (unit) length and opposite edges being parallel. We term these \textit{rectilinear} hexagons. We explored some other classes, but did not find a hexagonal tiling that gives a greater $d(k)$ than the best rectilinear one for any $k$ that we examined. 
However, our search was in no way exhaustive, so it remains very possible that such tilings exist.

\begin{figure}[!b]
\centering
\includegraphics[scale=0.2]{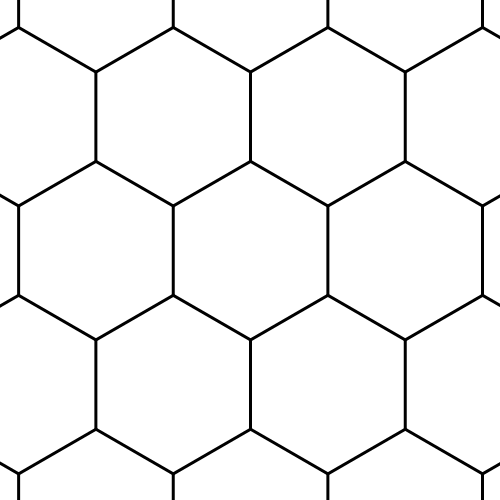}
\:\:\:
\includegraphics[scale=0.2]{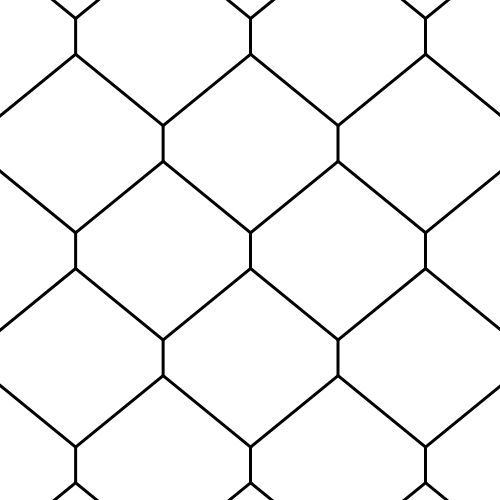}
\:\:\:
\includegraphics[scale=0.2]{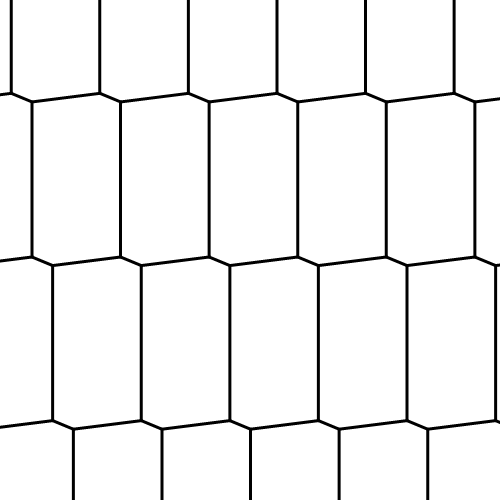}
\caption{Tilings of the plane with regular, semi-regular or rectilinear hexagons (respectively from left to right).}
\label{f2}
\end{figure}

Among this class of hexagon, there is a subset that we term \textit{semi-regular}: these have four edges of equal length, so they can be oriented with two vertices at $(0,\pm\frac{1}{2})$ and the others at $(\pm x,\pm y)$ with $x^2+y^2=\frac{1}{4}$. We decided to explore whether there exist values of $k$ for which semi-regular hexagons improve on regular hexagons, and similarly whether rectilinear hexagons ever beat semi-regular ones (see Fig.~\ref{f2}). It turns out that both such scenarios often occur.

\begin{figure}[!b]
\centering
\includegraphics[scale=0.24]{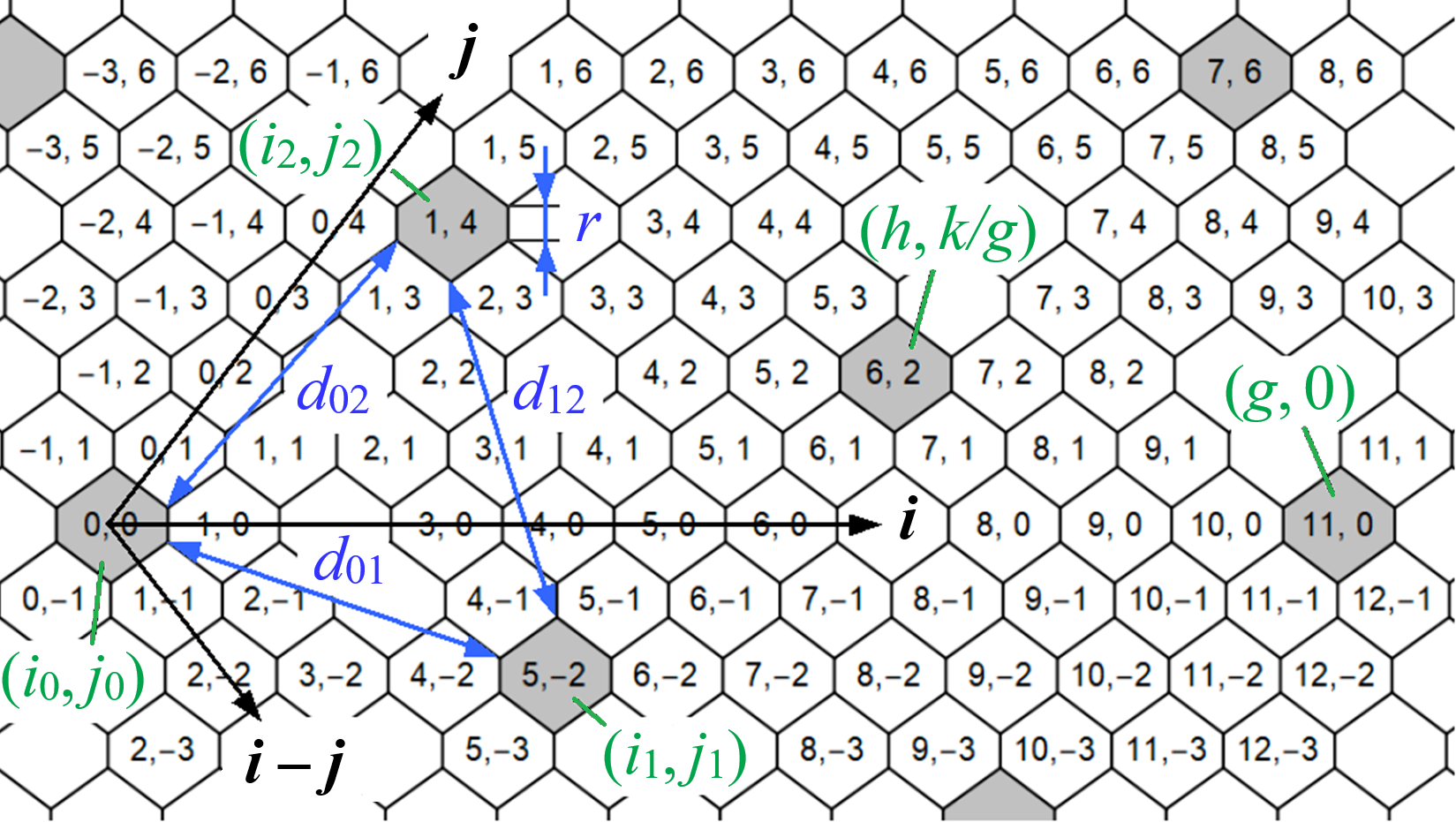} 
\caption{Notation for the tilings discussed in this study. In the example shown, $k=22$, $g=11$ and $h=6$.}
\label{f3}
\end{figure}

\paragraph{Tiling parameters.}
An example of tiling is shown in Fig.~\ref{f3}. Tiles of some (\textit{base}) colour are highlighted in grey. An oblique coordinate system $(i, j)$, where $i, j\in\mathbb{Z}$, is introduced for indexing tiles. For ease of visualization, one of the axes is oriented horizontally.

In our case, the specific tiling is determined by the position of two tiles of the base colour relative to the base tile, labeled $(0, 0)$. For calculations, it is convenient to use the tiles $(g, 0)$ and $(h, k/g)$. (We will define the parameters $g$ and $h$ below.) To present the results, it is more convenient to use the tiles $(i_1, j_1)$ and $(i_2, j_2)$, which together with the tile $(0, 0)$ form the triple with the smallest distances $\{d_{01}, d_{02}, d_{12}\}$ for the given tiling\footnote{
We used a similar construction in \cite{par} for the case of colouring the plane in several layers. But in that work we stipulated equality of distances in each triple of the nearest tiles, $d_{01}=d_{02}=d_{12}$. Perhaps some results obtained there can be improved by relaxing this restriction.},
$\min (d_{01}, d_{02}, d_{12})=d(k)$. The relation is fulfilled: $k=i_1 j_2-i_2 j_1$.

The shape of a tile is described by two parameters, for example, the length of two adjacent edges $\{r, s\}$, of which we are interested in the length $r$ of the vertical edge that completely defines the shape of a semi-regular hexagon.


\paragraph{Finding the shortest distance between hexagons.}

The colouring options for a given $k$ and a given class of hexagon can be defined in terms of only two parameters, $g$ and $h$, with $g$ a divisor of $k$ and $0 \le h < g$. Each horizontal row 
is coloured in a repeating sequence of $g$ colours, so for any $n$, row $k/g + n$ uses the same colours as row $n$, offset by $h$ hexagons. 

We are thus interested in the minimum distance between a point in hexagon $(0,0)$ and a point in hexagon $(ah+bg, ak/g)$, $a, b\in \mathbb{Z}$. It suffices to consider all tiles with coordinates $|\,i\,|+|\,j\,|<2\sqrt{k},\; j\ge 0$. We wish to identify the shape of hexagon (within the class being considered) and values of $g$ and $h$ that maximise that distance. For the classes we are considering, all of a hexagon’s angles must be obtuse, and each pair of opposite edges forms a rectangle. Thus, as with regular hexagons, $d(k)$ is either a multiple of the length of the relevant rectangle or a distance between corners of the two tiles.

\paragraph{Computation.}

Our work was performed using \texttt{Mathematica}, version 13. Briefly, we used numerical optimisation (the function \texttt{NMaxValue} or \texttt{FindMaximum}) to determine the best parameter values for a given $k$, and then, if the result improved upon any for a more restricted class of hexagon with the same or a smaller $k$ or for any class with a strictly smaller $k$, we derived exact expressions for $d(k)$ and the parameters using \texttt{Maximize}.

\section{Results}

\paragraph{Main results.}

Table~\ref{sult} compiles our findings, for $3 \le k \le 175$ (below we will assume this range of $k$ by default). Note that we show $d^2(k)$, not $d(k)$, because $d^2(k)$ is usually rational whereas $d(k)$ rarely is. A blank cell indicates that the best result obtained is no better than that for a more restricted class of hexagon with the same $k$ or for some class with a smaller $k$. However, for each $k$ we show the best result even when it is beaten by a smaller $k$; cases that do beat all smaller $k$ are denoted by bold-face in the last column.

The coordinates $(i_1, j_1)$, $(i_2, j_2)$ correspond to the tiling that provides the best obtained value of $d(k)$ for a given $k$. For semi-regular hexagons we also give a parameter $r$, which is the length of the vertical edges. Where possible, we impose additional restrictions on the parameters: tile 1 lies between the positive directions of the $(i-j)$ and $i$ axes, tile 2 lies between the $i$ and $j$ axes, and also the condition $d_{01}=d_{02}\le d_{12}$ is satisfied. In the case of $k=77$, these conditions contradict each other, and the distance $d_{01}(77)$ is the largest in the triplet.

The values for rectilinear $k=11,23,45$ and $187$ are roots of quartic equations. See the next subsection for more details, including the definition of the function $f(k)$ given as $d^2(k)$ for certain $k$.


\begin{table}[H]
{
\smallskip
\caption{Optimal $d(k)$ for regular, semi-regular and rectilinear hexagons.\label{sult}}
\vspace{1mm}

{
\centering
\footnotesize
\begin{tabular}[!t]{@{\,}c@{\;\,}|@{\;}c@{\;}c@{\;}c@{\;}c@{\;}| @{}P{0.10\linewidth}@{} | @{}P{0.28\linewidth}@{} | @{}P{0.305\linewidth}@{} | @{}r@{\,}}
\hline
\small{$k$} & \small{$i_1$} & \small{$j_1$} & \small{$i_2$} & \small{$j_2$} & \multicolumn{3}{c|}{\small{$d^2(k)$ [; $r$]}} & \multicolumn{1}{c}{\small{$d(k)$}} \\
\cline{6-8}
 &&&&& \small{regular} & \small{semi-regular} & \small{rectilinear} & \multicolumn{1}{c}{\small{approx.}} \\
\hline \hline
3   &  2 & -1 & 1 &  1 & $(1/2)^2$ &  &  &\m 0.500000 \\ \hline
4   &  2 &  0 & 0 &  2 & 3/4 &  &  &\m 0.866025 \\ \hline
5   &  3 & -1 & 2 &  1 &  & $(5/6)^2$; 5/6 &  & 0.833333 \\ \hline
6   &  3 &  0 & 0 &  2 &  & $\frac{17 + \sqrt{33}}{24}$; $\frac{\sqrt{33} - 3}{12}$ & $f(6)$ &\m 0.992076 \\ \hline
7   &  3 & -1 & 1 &  2 & 7/4 &  &  &\m 1.322876 \\ \hline
8   &  3 & -2 & 1 &  2 &  & $(7/5)^2$; 1/5 &  &\m 1.400000 \\ \hline
9   &  3 &  0 & 0 &  3 & 3 &  &  &\m 1.732051 \\ \hline
10  &  4 & -1 & 2 &  2 &  & 11/4; 3/4 &  & 1.658312 \\ \hline
11  &  3 & -2 & 1 &  3 &  & 11/4; 3/4 & $1.674155^2$ & 1.674155 \\ \hline
12  &  4 & -2 & 2 &  2 & $2^2$ &  &  &\m 2.000000 \\ \hline
13  &  4 & -1 & 1 &  3 & 19/4 &  &  &\m 2.179449 \\ \hline
14  &  4 & -2 & 1 &  3 &  & 19/4; 1/4 &  & 2.179449 \\ \hline
15  &  4 & -3 & 1 &  3 &  & 153/32; 3/8 &  &\m 2.186607 \\ \hline
16  &  4 &  0 & 0 &  4 & 27/4 &  &  &\m 2.598076 \\ \hline
17  &  5 & -1 & 2 &  3 &  & 627/100; 7/10 & 709371/111724 & 2.519785 \\ \hline
18  &  4 & -3 & 2 &  3 &  & 6; 2/3 & $20-\sqrt{192}$ & 2.478627 \\ \hline
19  &  5 & -2 & 2 &  3 & 31/4 &  &  &\m 2.783882 \\ \hline
20  &  5 &  0 & 0 &  4 &  & $\frac{5931 +405 \sqrt{129}}{1352}$; $\frac{3 \sqrt{129} - 15}{52}$ & $f(20)$ &\m 2.862821 \\ \hline
21  &  5 & -1 & 1 &  4 & 37/4 &  &  &\m 3.041381 \\ \hline
22  &  5 & -2 & 1 &  4 &  & $3^2$; 1/3 & 193496/21275 & 3.015791 \\ \hline
23  &  5 & -3 & 1 &  4 &  & 103/12; 1/6 & $2.942985^2$ & 2.942985 \\ \hline
24  &  6 & -2 & 3 &  3 &  & 1216/121; 29/44 & 869/86 &\m 3.178781 \\ \hline
25  &  5 &  0 & 0 &  5 & 12 &  &  &\m 3.464102 \\ \hline
26  &  6 & -1 & 2 &  4 &  & 203/18; 2/3 & 2150/187 & 3.390771 \\ \hline
27  &  6 & -3 & 3 &  3 & $(7/2)^2$ &  &  &\m 3.500000 \\ \hline
28  &  6 & -2 & 2 &  4 & 13 &  &  &\m 3.605551 \\ \hline
29  &  7 & -2 & 4 &  3 &  & 44719/3844; 47/62 & 481078403/41028908 & 3.424230 \\ \hline
30  &  6 &  0 & 0 &  5 &  & $\frac{1072 + 128 \sqrt{51}}{147}$; $\frac{2\sqrt{51} - 6}{21}$ & $f(30)$ &\m 3.762292 \\ \hline
31  &  6 & -1 & 1 &  5 & 61/4 &  &  &\m 3.905125 \\ \hline
32  &  6 & -1 & 2 &  5 &  & 59/4; 3/8 & 2307733/153928 & 3.871988 \\ \hline
33  &  7 & -3 & 4 &  3 &  & $(227/58)^2$; 37/58 &  &\m 3.913793 \\ \hline
34  &  7 & -2 & 3 &  4 &  & 1921/121; 7/11 & 170126/10611 &\m 4.004121 \\ \hline
35  &  7 &  0 & 1 &  5 &  \multicolumn{2}{r@{\;}|}{$\frac{17559 + 135 \sqrt{12257}}{2048}; \frac{\sqrt{12257} - 15}{128}$} &  & 3.983917 \\ \hline
36  &  6 &  0 & 0 &  6 & 75/4 &  &  &\m 4.330127 \\ \hline
37  &  7 & -3 & 3 &  4 & 73/4 &  &  & 4.272002 \\ \hline
38  &  6 & -1 & 2 &  6 &  & 17; 3/5 & 343577/19520 & 4.195388 \\ \hline
39  &  7 & -2 & 2 &  5 & 79/4 &  &  &\m 4.444097 \\ \hline
40  &  7 & -4 & 3 &  4 &  & $(137/31)^2$; 11/31 &  & 4.419355 \\ \hline
41  &  7 & -2 & 3 &  5 &  & 637/32; 3/8 & 19168435/946444 &\m 4.500346 \\ \hline
42  &  7 &  0 & 0 &  6 &  \multicolumn{2}{r@{\;}|}{$\frac{84075 + 13125 \sqrt{33}}{7688}$; $\frac{15 \sqrt{33} - 35}{124}$} & $f(42)$ &\m 4.651075 \\ \hline
43  &  7 & -1 & 1 &  6 & 91/4 &  &  &\m 4.769696 \\ \hline
44  &  7 & -1 & 2 &  6 &  & 22; 2/5 & 4700717/209612 & 4.735589 \\ \hline
45  &  7 & -3 & 1 &  6 &  & 423/20; 3/10 & $4.607942^2$ & 4.607942 \\ \hline
\end{tabular}
}
}
\end{table}

\begin{table}[H]
{
\smallskip
\centering
Continuation of Table \ref{sult}. 
\vspace{2mm}

{
\centering
\footnotesize
\begin{tabular}[!t]{@{\,}c@{\;\,}|@{\;}c@{\;}c@{\;}c@{\;}c@{\;}| @{}P{0.10\linewidth}@{} | @{}P{0.28\linewidth}@{} | @{}P{0.305\linewidth}@{} | @{}r@{\,}}
\hline
\small{$k$} & \small{$i_1$} & \small{$j_1$} & \small{$i_2$} & \small{$j_2$} & \multicolumn{3}{c|}{\small{$d^2(k)$ [; $r$]}} & \multicolumn{1}{c}{\small{$d(k)$}} \\
\cline{6-8}
 &&&&& \small{regular} & \small{semi-regular} & \small{rectilinear} & \multicolumn{1}{c}{\small{approx.}} \\
\hline \hline
46  &  8 & -2 & 3 &  5 &  & 49051/2116; 57/92 & 73151/3116 &\m 4.845197 \\ \hline
47  &  7 & -3 & 4 &  5 &  & 3224/147; 13/21 & 5560412771/243356228 & 4.780048 \\ \hline
48  &  8 & -4 & 4 &  4 & $5^2$ &  &  &\m 5.000000 \\ \hline
49  &  7 &  0 & 0 &  7 & 27 &  &  &\m 5.196152 \\ \hline
50  &  8 & -1 & 2 &  6 &  & 825/32; 5/8 & 2143125/81224 & 5.136669 \\ \hline
51  &  7 & -1 & 2 &  7 &  & 6016/243; 19/27 & 32914863743/1281311828 & 5.068374 \\ \hline
52  &  8 & -2 & 2 &  6 & 28 &  &  &\m 5.291503 \\ \hline
53  &  9 & -2 & 4 &  5 &  & 358639/13924; 83/118 & 21713005483/826893628 & 5.124307 \\ \hline
54  &  8 & -2 & 3 &  6 &  & 335/12; 5/12 & 118349/4160 &\m 5.333787 \\ \hline
55  &  9 & -1 & 1 &  6 &  & $(449/86)^2$; 21/86 & 15510042323/562475260 & 5.251154 \\ \hline
56  &  8 &  0 & 0 &  7 &  & $(277/51)^2$; 31/51 & $f(56)$ &\m 5.533452 \\ \hline
57  &  8 & -1 & 1 &  7 & 127/4 &  &  &\m 5.634714 \\ \hline
58  &  8 & -1 & 2 &  7 &  & 123/4; 5/12 & 7191956/229175 & 5.601959 \\ \hline
59  &  8 & -1 & 3 &  7 &  & 89/3; 1/3 & 742918955/24763644 & 5.477261 \\ \hline
60  &  9 & -2 & 3 &  6 &  & 118969/3721; 37/61 & 4775023/147232 &\m 5.694907 \\ \hline
61  &  9 & -4 & 4 &  5 & 133/4 &  &  &\m 5.766281 \\ \hline
62  &  8 & -2 & 3 &  7 &  & 49987/1600; 47/80 & 2563175/79001 & 5.696037 \\ \hline
63  &  9 & -3 & 3 &  6 & 139/4 &  &  &\m 5.894913 \\ \hline
64  &  8 &  0 & 0 &  8 & 147/4 &  &  &\m 6.062178 \\ \hline
65  &  9 & -1 & 2 &  7 &  & 11431/324; 11/18 & 42769993/1184612 & 6.008714 \\ \hline
66  &  9 & -3 & 4 &  6 &  & 23987/675; 17/45 & 1193079661/33188989 & 5.995670 \\ \hline
67  &  9 & -2 & 2 &  7 & 151/4 &  &  &\m 6.144103 \\ \hline
68  & 10 & -2 & 4 &  6 &  & 12649/361; 13/19 & 3193957/89125 & 5.986387 \\ \hline
69  &  9 & -2 & 3 &  7 &  & 262/7; 3/7 & 6272261/164372 &\m 6.177292 \\ \hline
70  & 10 & -4 & 5 &  5 &  & 44339/1156; 81/136 & 414848/10783 &\m 6.202613 \\ \hline
71  & 10 & -1 & 1 &  7 &  & 1035/28; 5/14 & 3274003629/86168876 & 6.164025 \\ \hline
72  &  9 &  0 & 0 &  8 &  \multicolumn{2}{r@{\;}|}{$\frac{177233 + 7203 \sqrt{537}}{8664}$; $\frac{7 \sqrt{537} - 63}{228}$} & $f(72)$ &\m 6.411744 \\ \hline
73  &  9 & -1 & 1 &  8 & $(13/2)^2$ &  &  &\m 6.500000 \\ \hline
74  &  9 & -1 & 2 &  8 &  & 41; 3/7 & 9026779/215683 & 6.469317 \\ \hline
75  & 10 & -5 & 5 &  5 & $(13/2)^2$ &  &  & 6.500000 \\ \hline
76  & 10 & -4 & 4 &  6 & 43 &  &  &\m 6.557439 \\ \hline
77  & 11 &  0 & 2 &  7 &  & 1215/32; 5/8 &  & 6.161879 \\ \hline
78  &  9 & -2 & 3 &  8 &  & 21844/529; 53/92 & 4833828787/112826683 & 6.545453 \\ \hline
79  & 10 & -3 & 3 &  7 & 181/4 &  &  &\m 6.726812 \\ \hline
80  & 10 & -4 & 5 &  6 &  \multicolumn{2}{r@{\;}|}{$\frac{63261 + 972 \sqrt{3770}}{2738}$; $\frac{\sqrt{3770} - 12}{74}$} & 241351/5371 & 6.703428 \\ \hline
81  &  9 &  0 & 0 &  9 & 48 &  &  &\m 6.928203 \\ \hline
82  & 10 & -1 & 2 &  8 &  & 1157/25; 3/5 & 1184696/25029 & 6.879893 \\ \hline
83  &  9 & -1 & 2 &  9 &  & 957559/21316; 47/146 & \scriptsize{975386934409/20985431116} & 6.817569 \\ \hline
84  & 10 & -2 & 2 &  8 & $7^2$ &  &  &\m 7.000000 \\ \hline
85  & 11 & -5 & 6 &  5 &  & $(1097/158)^2$; 93/158 &  & 6.943038 \\ \hline
86  & 10 & -2 & 3 &  8 &  & 40953/841; 17/29 & 30121/610 &\m 7.026997 \\ \hline
87  &  9 & -2 & 3 &  9 &  & 638661/13924; 79/118 & 1080643/23188 & 6.826679 \\ \hline
88  & 10 & -2 & 4 &  8 &  & 191/4; 3/8 & 84065/1729 & 6.972847 \\ \hline
89  & 11 & -3 & 4 &  7 &  & 194709/3844; 109/186 & \scriptsize{511258004631/9991926364} &\m 7.153119 \\ \hline
\end{tabular}
}
}
\end{table}

\begin{table}[H]
{
\smallskip
\centering
Continuation of Table \ref{sult}.
\vspace{2mm}

{
\centering
\footnotesize
\begin{tabular}[!t]{@{}c@{\;}|@{\;}c@{\;}c@{\;}c@{\;}c@{\;}| @{}P{0.10\linewidth}@{} | @{}P{0.28\linewidth}@{} | @{}P{0.305\linewidth}@{} | @{}r@{\,}}
\hline
\small{$k$} & \small{$i_1$} & \small{$j_1$} & \small{$i_2$} & \small{$j_2$} & \multicolumn{3}{c|}{\small{$d^2(k)$ [; $r$]}} & \multicolumn{1}{c}{\small{$d(k)$}} \\
\cline{6-8}
 &&&&& \small{regular} & \small{semi-regular} & \small{rectilinear} & \multicolumn{1}{c}{\small{approx.}} \\
\hline \hline
90  & 10 &  0 & 0 &  9 &  \multicolumn{2}{r@{\;}|}{$\frac{140352 + 30720 \sqrt{19}}{5329}$; $\frac{12 \sqrt{19} - 20}{73}$} & $f(90)$ &\m 7.287290 \\ \hline
91  & 10 & -1 & 1 &  9 & 217/4 &  &  &\m 7.365460 \\ \hline
92  & 10 & -1 & 2 &  9 &  & 211/4; 7/16 & 4639657267/86189275 & 7.336963 \\ \hline
93  & 11 & -4 & 4 &  7 & 217/4 &  &  & 7.365460 \\ \hline
94  & 11 & -2 & 3 &  8 &  & $(713/97)^2$; 57/97 & 666148999/12139951 &\m 7.407595 \\ \hline
95  & 12 & -5 & 7 &  5 &  & $(277/38)^2$; 25/38 &  & 7.289474 \\ \hline
96  & 11 & -6 & 5 &  6 &  & $(67/9)^2$; 11/27 &  &\m 7.444444 \\ \hline
97  & 11 & -3 & 3 &  8 & 229/4 &  &  &\m 7.566373 \\ \hline
98  & 12 & -1 & 2 &  8 &  & 115625/2178; 8/33 & 613512713/11338873 & 7.355748 \\ \hline
99  & 11 &  0 & 1 &  9 &  \multicolumn{2}{r@{\;}|}{$\frac{498775+1575 \sqrt{91553}}{16928}$; $\frac{\sqrt{91553} - 63}{368}$} &  &\m 7.590561 \\ \hline
100 & 10 &  0 & 0 & 10 & 243/4 &  &  &\m 7.794229 \\ \hline
101 & 11 & -1 & 2 &  9 &  & 28449/484; 13/22 & 12663270675/210818572 & 7.750300 \\ \hline
102 & 12 & -5 & 6 &  6 &  & 29651/500; 29/50 & 208565/3509 & 7.709550 \\ \hline
103 & 11 & -2 & 2 &  9 & 247/4 &  &  &\m 7.858117 \\ \hline
104 & 12 & -4 & 5 &  7 &  & 176485/2916; 125/216 & 4142419871/67995631 & 7.805245 \\ \hline
105 & 11 & -2 & 3 &  9 &  & 1097/18; 4/9 & 447216991/7200892 &\m 7.880722 \\ \hline
106 & 10 & -2 & 3 & 10 &  & $(84/11)^2$; 29/44 & 33151587/560347 & 7.691723 \\ \hline
107 & 11 & -5 & 6 &  7 &  & 2161/36; 7/18 & $\frac{12486512240719}{202886322964}$ & 7.845022 \\ \hline
108 & 12 & -6 & 6 &  6 & $8^2$ &  &  &\m 8.000000 \\ \hline
109 & 12 & -5 & 5 &  7 & 259/4 &  &  &\m 8.046738 \\ \hline
110 & 11 &  0 & 0 & 10 &  \multicolumn{2}{r@{\;}|}{$\frac{2184003+72171 \sqrt{849}}{66248}$; $\frac{9 \sqrt{849} - 99}{364}$} & $f(110)$ &\m 8.160909 \\ \hline
111 & 11 & -1 & 1 & 10 & 271/4 &  &  &\m 8.231039 \\ \hline
112 & 11 & -1 & 2 & 10 & 67 &  & 1034089141/15361771 & 8.204618 \\ \hline
113 & 11 & -1 & 3 & 10 &  & 451319/6962; 38/59 & \scriptsize{387649933983/5906484172} & 8.101312 \\ \hline
114 & 12 & -2 & 3 &  9 &  & 937099/13924; 137/236 & 283698767/4150532 &\m 8.267550 \\ \hline
115 & 13 & -4 & 6 &  7 &  & \scriptsize{4244941/64516; 163/254} & $\frac{12272144911631}{184524424676}$ & 8.155175 \\ \hline
116 & 12 & -4 & 5 &  8 &  & 8723/128; 13/32 & 2045764331/29695775 &\m 8.300045 \\ \hline
117 & 12 & -3 & 3 &  9 & 283/4 &  &  &\m 8.411302 \\ \hline
118 & 11 & -2 & 4 & 10 &  & 451319/6962; 38/59 & 22520512223/336563063 & 8.180049 \\ \hline
119 & 13 & -3 & 5 &  8 &  & $(1643/198)^2$; 65/198 & $\frac{17514883072583}{252372028148}$ & 8.330729 \\ \hline
120 & 12 & -3 & 4 &  9 &  & $\frac{28677 + 3630 \sqrt{58}}{784}$; $\frac{3 \sqrt{58} - 5}{28}$ & \scriptsize{102121645769/1420956281} &\m 8.477515 \\ \hline
121 & 11 &  0 & 0 & 11 & 75 &  &  &\m 8.660254 \\ \hline
122 & 12 & -1 & 2 & 10 &  & 655/9; 7/12 & 10143575/136512 & 8.620056 \\ \hline
123 & 11 & -1 & 2 & 11 &  & 392863/5476; 79/222 & \scriptsize{129728247871/1768473940} & 8.564814 \\ \hline
124 & 12 & -2 & 2 & 10 & 76 &  &  &\m 8.717798 \\ \hline
125 & 13 & -2 & 4 &  9 &  & \scriptsize{5551687/77284; 179/278} & \scriptsize{976273547119/13237256716} & 8.587896 \\ \hline
126 & 12 & -2 & 3 & 10 &  & 1499/20; 9/20 & 19145005/250792 &\m 8.737172 \\ \hline
127 & 13 & -6 & 6 &  7 & 307/4 &  &  &\m 8.760708 \\ \hline
128 & 12 & -4 & 5 &  9 &  & 369/5; 2/5 & 795830295/10506727 & 8.703151 \\ \hline
129 & 13 & -5 & 5 &  8 & 313/4 &  &  &\m 8.845903 \\ \hline
130 & 13 &  0 & 1 & 10 &  \multicolumn{2}{r@{\;}|}{$\frac{3182571+312741 \sqrt{97}}{84872}$; $\frac{27 \sqrt{97} - 143}{412}$} & $\frac{18088 + 16 \sqrt{1192555}}{469}$ & 8.707601 \\ \hline
131 & 13 & -1 & 1 & 10 &  & 3094091/40804; 81/202 & \scriptsize{192224305737/2452589252} &\m 8.853026 \\ \hline
132 & 12 &  0 & 0 & 11 &  \multicolumn{2}{r@{\;}|}{$\frac{165700 + 10000 \sqrt{258}}{4107}$; $\frac{5 \sqrt{258} - 30}{111}$} & $f(132)$ &\m 9.033127 \\ \hline
\end{tabular}
}
}
\end{table}

\begin{table}[H]
{
\smallskip
\centering
Continuation of Table \ref{sult}. 
\vspace{2mm}

{
\centering
\footnotesize
\begin{tabular}[!t]{@{}c@{\;}|@{\;}c@{\;}c@{\;}c@{\;}c@{\;}| @{}P{0.10\linewidth}@{} | @{}P{0.28\linewidth}@{} | @{}P{0.305\linewidth}@{} | @{}r@{\,}}
\hline
\small{$k$} & \small{$i_1$} & \small{$j_1$} & \small{$i_2$} & \small{$j_2$} & \multicolumn{3}{c|}{\small{$d^2(k)$ [; $r$]}} & \multicolumn{1}{c}{\small{$d(k)$}} \\
\cline{6-8}
 &&&&& \small{regular} & \small{semi-regular} & \small{rectilinear} & \multicolumn{1}{c}{\small{approx.}} \\
\hline \hline
133 & 12 & -1 & 1 & 11 & 331/4 &  &  &\m 9.096703 \\ \hline
134 & 12 & -1 & 2 & 11 &  & 323/4; 9/20 & 6558879037/79690552 & 9.072174 \\ \hline
135 & 12 & -1 & 3 & 11 &  & 787/10; 2/5 & $\frac{10272242758901}{127501682004}$ & 8.975831 \\ \hline
136 & 13 & -2 & 3 & 10 &  & 181289/2209; 27/47 & 20689231/248261 &\m 9.128889 \\ \hline
137 & 13 & -4 & 5 &  9 &  & 2598661/31684; 73/178 & $\frac{74222844365281}{892932274564}$ & 9.117159 \\ \hline
138 & 12 & -2 & 3 & 11 &  & 27389/338; 19/52 & 1122833/13505 & 9.118225 \\ \hline
139 & 13 & -3 & 3 & 10 & 343/4 &  &  &\m 9.260130 \\ \hline
140 & 14 & -6 & 7 &  7 &  & 134518/1587; 157/276 & 80581/949 & 9.214743 \\ \hline
141 & 14 & -3 & 5 &  9 &  & 929347/11236; 67/106 & \scriptsize{196014091823/2324718468} & 9.182448 \\ \hline
142 & 13 & -3 & 4 & 10 &  & 117688/1369; 21/37 & 15150070069/174437335 &\m 9.319392 \\ \hline
143 & 13 &  0 & 1 & 11 &  \multicolumn{2}{r@{\;}|}{$\frac{399015 + 4455 \sqrt{7553}}{8978}$; $\frac{5 \sqrt{7553}-99}{536}$} &  &\m 9.357805 \\ \hline
144 & 12 &  0 & 0 & 12 & 363/4 &  &  &\m 9.526279 \\ \hline
145 & 13 & -1 & 2 & 11 &  & 59675/676; 15/26 & 3279780625/36423236 & 9.489277 \\ \hline
146 & 14 & -4 & 5 &  9 &  & 119853/1352; 59/104 & 16400637/183488 & 9.454238 \\ \hline
147 & 13 & -2 & 2 & 11 & 367/4 &  &  &\m 9.578622 \\ \hline
148 & 14 & -6 & 6 &  8 & 91 &  &  & 9.539392 \\ \hline
149 & 13 & -2 & 3 & 11 &  & 995/11; 5/11 & 1164694451/12649484 &\m 9.595544 \\ \hline
150 & 13 & -4 & 5 & 10 &  & 112529/1280; 103/160 & \scriptsize{286692851539/3160973281} & 9.523531 \\ \hline
151 & 14 & -5 & 5 &  9 & 373/4 &  &  &\m 9.656604 \\ \hline
152 & 14 & -3 & 4 & 10 &  & 405821/4374; 46/81 & 7532633043/80257313 &\m 9.687932 \\ \hline
153 & 15 & -6 & 8 &  7 &  & \scriptsize{1030243/11236; 199/318} & $\frac{894213901557691}{9724451350300}$ & 9.589328 \\ \hline
154 & 14 & -6 & 7 &  8 &  & 771407/8192; 109/256 & 1311592/13897 &\m 9.714912 \\ \hline
155 & 14 & -1 & 1 & 11 &  & 133856/1445; 53/85 & $\frac{15557710498703}{164015315380}$ &\m 9.739365 \\ \hline
156 & 13 &  0 & 0 & 12 & 97 &  & $f(156)$ &\m 9.904296 \\ \hline
157 & 13 & -1 & 1 & 12 & 397/4 &  &  &\m 9.962429 \\ \hline
158 & 13 & -1 & 2 & 12 &  & 97; 5/11 & 75425702171/763451750 & 9.939599 \\ \hline
159 & 13 & -1 & 3 & 12 &  \multicolumn{2}{r@{\;}|}{\scriptsize{12176361/128164; 223/358}} & \scriptsize{184097879615/1897621332} & 9.849623 \\ \hline
160 & 14 & -2 & 3 & 11 &  & 677416/6889; 189/332 & 9726226252/97433683 &\m 9.991199 \\ \hline
161 & 15 & -7 & 8 &  7 &  & \scriptsize{9284209/93636; 173/306} &  & 9.957516 \\ \hline
162 & 15 & -6 & 7 &  8 &  & 301273/3025; 31/55 & \scriptsize{195650535281/1958502269} &\m 9.994901 \\ \hline
163 & 14 & -3 & 3 & 11 & 409/4 &  &  &\m 10.11187 \\ \hline
164 & 14 & -4 & 6 & 10 &  & 292/3; 5/12 & 267614/2691 & 9.972351 \\ \hline
165 & 15 & -5 & 6 &  9 &  & \scriptsize{2487641/24500; 197/350} & $\frac{33252082489597}{325566557812}$ & 10.10624 \\ \hline
166 & 14 & -3 & 4 & 11 &  & 813/8; 11/24 & 544579139/5269928 &\m 10.16549 \\ \hline
167 & 15 & -1 & 2 & 11 &  & 6025357/60516; 77/246 & 7169589821/70929460 & 10.05388 \\ \hline
168 & 14 &  0 & 1 & 12 &  \multicolumn{2}{r@{\;}|}{$\frac{662311+5070\sqrt{16226}}{12482}$; $\frac{\sqrt{16226}-30}{158}$} &  &\m 10.23727 \\ \hline
169 & 13 &  0 & 0 & 13 & 108 &  &  &\m 10.39230 \\ \hline
170 & 14 & -1 & 2 & 12 &  & 10317/98; 4/7 & 1618378002/15084233 & 10.35806 \\ \hline
171 & 15 & -6 & 6 &  9 & 427/4 &  &  & 10.33199 \\ \hline
172 & 14 & -2 & 2 & 12 & 109 &  &  &\m 10.44031 \\ \hline
173 & 15 & -2 & 4 & 11 &  \multicolumn{2}{r@{\;}|}{\scriptsize{15153631/145924; 239/382}} & $\frac{37802632989571}{354544269100}$ & 10.32585 \\ \hline
174 & 14 & -2 & 3 & 12 &  & 2579/24; 11/24 & 1893749/17324 &\m 10.45531 \\ \hline
175 & 15 & -5 & 5 & 10 & 439/4 &  &  &\m 10.47616 \\ \hline
\end{tabular}
}
}
\end{table}


\begin{figure}[H]
\centering
{
\centering
\begin{tabular}{@{}c@{\;}cc@{\;}c@{}}
 3 & \includegraphics[scale=0.25]{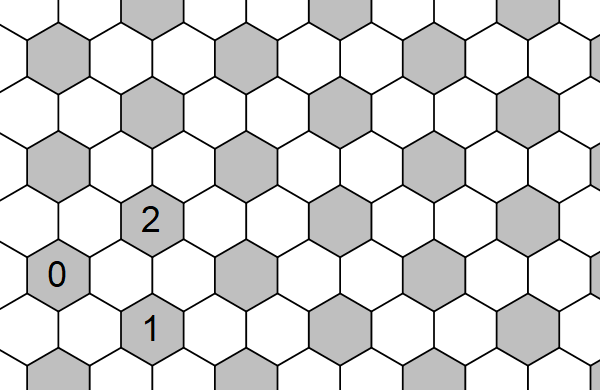}  & 7 & \includegraphics[scale=0.25]{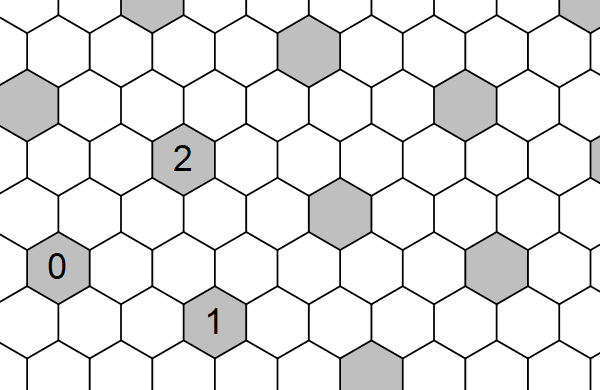}  \\ [2mm]
 4 & \includegraphics[scale=0.25]{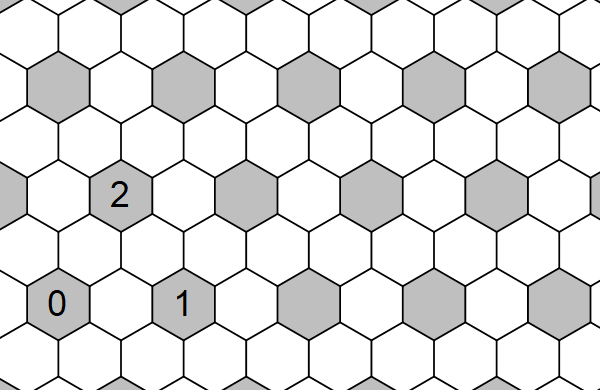}  & 8 & \includegraphics[scale=0.25]{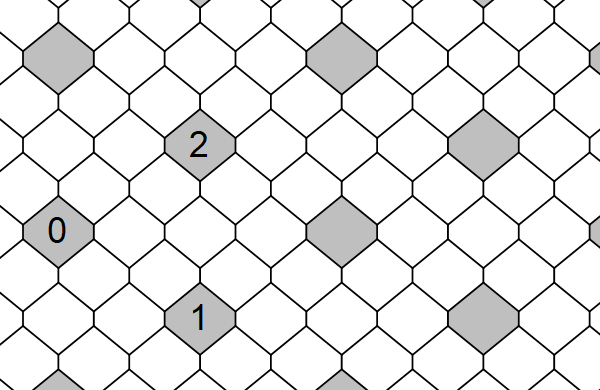}  \\ [2mm]
 5 & \includegraphics[scale=0.25]{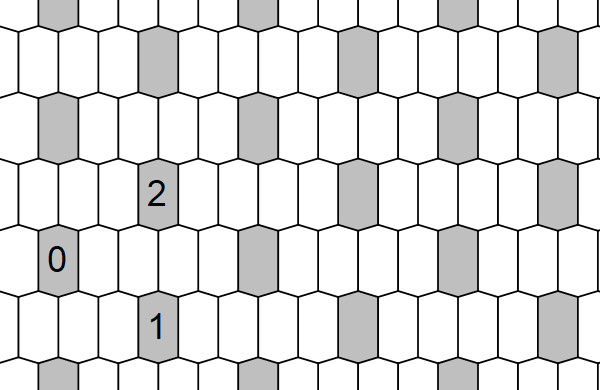}  & 9 & \includegraphics[scale=0.25]{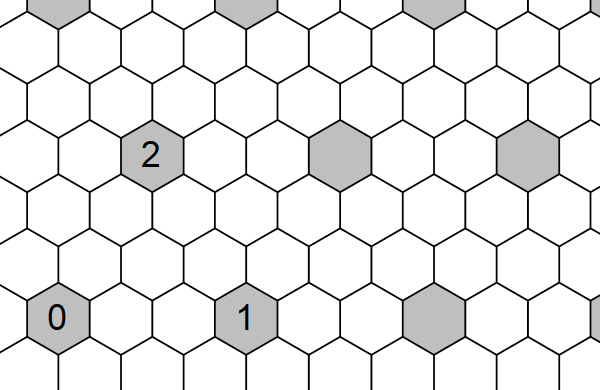}  \\ [2mm]
 6 & \includegraphics[scale=0.25]{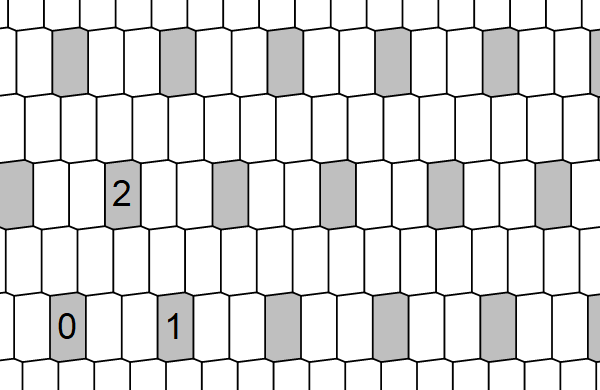}  & 10& \includegraphics[scale=0.25]{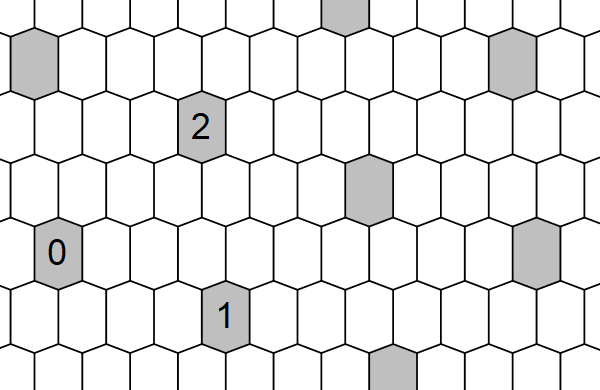} \\ [2mm]
*6 & \includegraphics[scale=0.25]{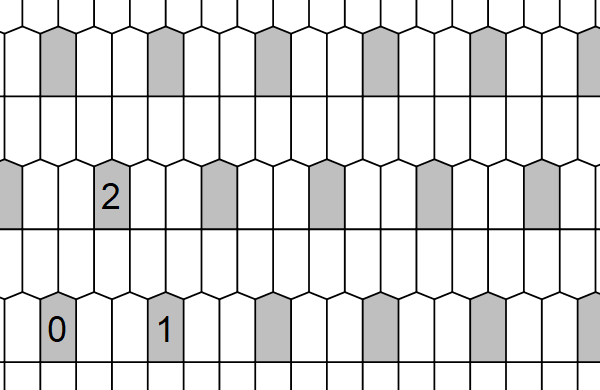}  & 11& \includegraphics[scale=0.25]{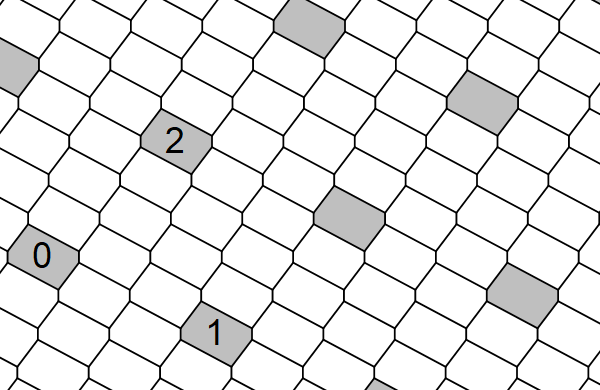} \\
\end{tabular} \par
}
\caption{Optimal hexagonal (* and pentagonal) tilings for $3 \leq k \leq 11$. }
\label{hex1}
\end{figure}

\begin{figure}[H]
\centering
{
\centering
\begin{tabular}{@{}c@{\;}cc@{\;}c@{}}
 12& \includegraphics[scale=0.25]{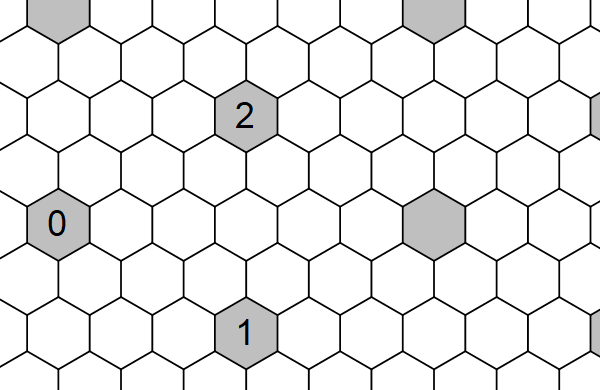} & 17& \includegraphics[scale=0.25]{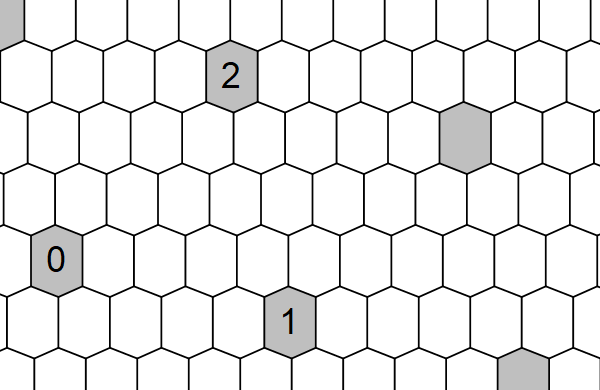} \\ [2mm]
 13& \includegraphics[scale=0.25]{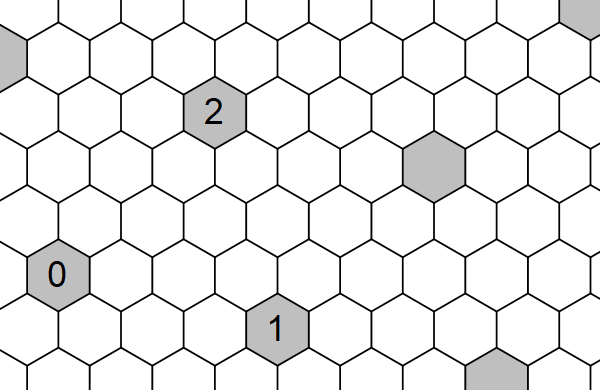} & 18& \includegraphics[scale=0.25]{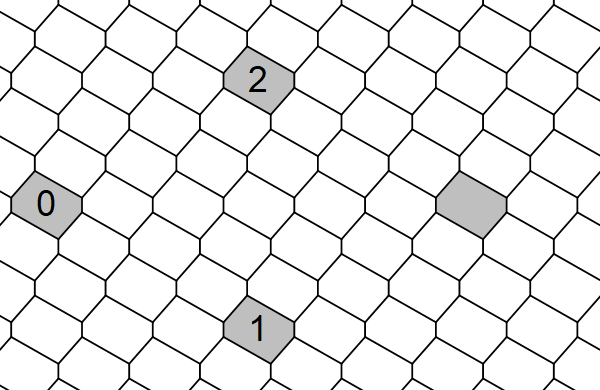} \\ [2mm]
 14& \includegraphics[scale=0.25]{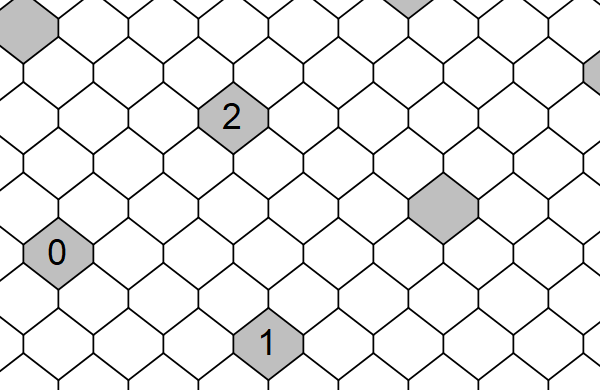} & 20& \includegraphics[scale=0.25]{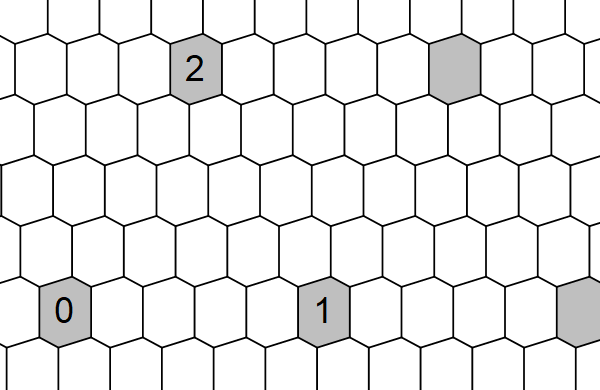} \\ [2mm]
 15& \includegraphics[scale=0.25]{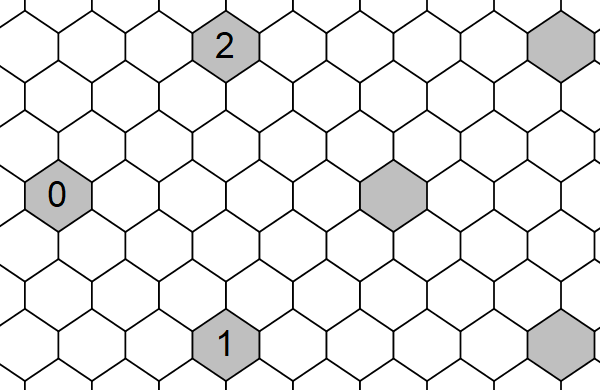} & 22& \includegraphics[scale=0.25]{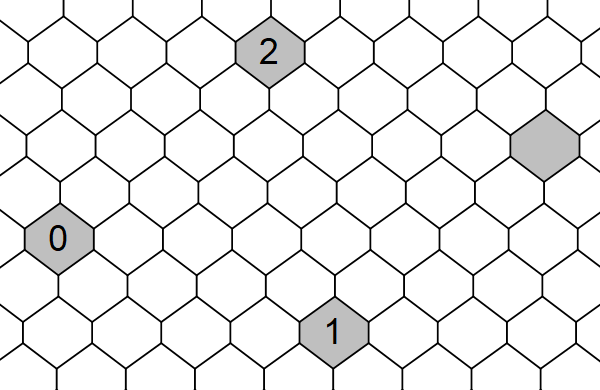} \\ [2mm]
 16& \includegraphics[scale=0.25]{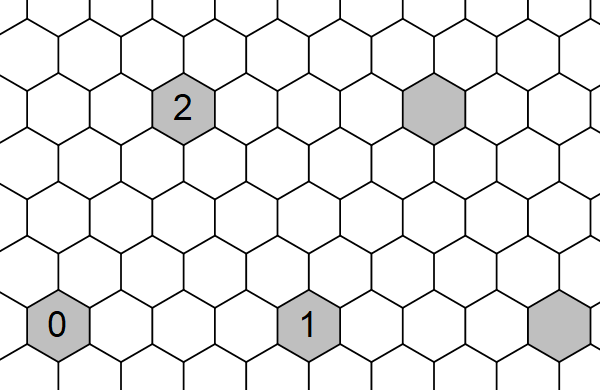} & 23& \includegraphics[scale=0.25]{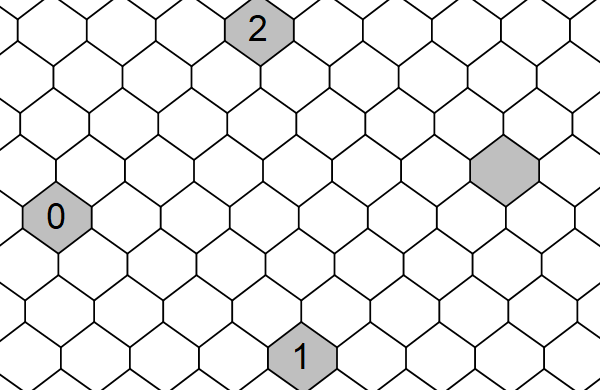} \\
\end{tabular} \par
}
\caption{Optimal hexagonal tilings for $12 \leq k \leq 23$. }
\label{hex2}
\end{figure}

Tilings with the maximum $d(k)$ that we obtained for various small $k$ are shown in Figures~\ref{hex1} and \ref{hex2}. All tilings are given in the same scale. Note that the hexagons with $k=17,22$ and $23$ are not precisely semi-regular.

\paragraph{Classification.}

Other than the L\"{o}schian class $k = 3n^2$ ($n \ge 1$), for which $d(k) = \frac{3n}{2} - 1$ in the regular setting (except for $k=147$, where another colouring of regular hexagons wins), $d(k)$ is not often rational; the only examples we found are regular $k=84$ and semi-regular $k=5$, 8, 22, 33, 40, 55, 56, 85, 94, 95, 96, 106 and 119 (some of which are beaten by rectilinear examples). 

However, $d^2(k)$ is usually rational.
The obtained values 
can be classified according to the degree of the polynomial of which $d^2(k)$ is a root:

i) When $d^2(k)$ is rational, in the semi-regular setting the denominator is always a square or a small multiple of a square, and is seldom much greater than $k^2$. (Also, the denominator of $r$ does not exceed $3k$.) In the rectilinear setting, however, the denominator is often far larger, reaching the tens of millions even for $k$ as small as $29$.

ii)	In all cases where the semi-regular $d^2(k)$ is irrational, it has the form $\mathbb{Q}[\sqrt{a}]$ for some integer $a$. All such solutions appear at $k=n(n+m)$, $0<m<<n$. In the rectilinear setting, there are a few examples where the best $d^2(k)$ has the quadratic form. In the range of $k$ that we explored, the only such examples are $k=18$ and $130$ in which there is a rectilinear tiling that beats the semi-regular one, and $k=35, 99$ and $143$ in which the quadratic-form semi-regular $d(k)$ is unbeaten. 

iii) The largest irrational class is $k=n(n+1)$ for $n \ge 2$, with the exception of the L\"{o}schian case $k=12$. Here, 
the rectilinear $d^2(k)$, 
denoted in the table by $f(k)$, is the largest real root\footnote{The other two roots of $f(k)$ tend to $-1/3$ and $+1$ as $k$ increases.} 
of the cubic polynomial 
$a_3x^3 + a_2x^2 + a_1x + a_0 = 0$, where $p=n(n-1)=k+1-\sqrt{4k+1}$, $a_3=4p(p^2+3p+1)$, $a_2=-3p^4-8p^3+2p^2+4p+1$, $a_1=2p^2(p^2-2p-1)$, $a_0=p^4$.
We do not know why. 

iv)	Only four values that we checked, namely $k=11, 23, 45$ and $187$, have an optimal $x=d^2(k)$ in the rectilinear setting that is a root of a polynomial of degree greater than 3, in all cases quartic. For completeness we list those polynomials here (take the smaller of the two real roots):

$25600 x^4 - 1459616 x^3 + 8840377 x^2 - 18735876 x + 13623552$

$97344 x^4 - 14493200 x^3 + 287680857 x^2 - 2041299600 x + 4968218624$

$548800 x^4 - 70030800 x^3 + 2721532527 x^2 - 44612184348 x + 279110188800$

\abovedisplayskip=-10pt

\begin{multline*}
\;\;391936098304 x^4 - 214522652834752 x^3 + 43720132398169569 x^2- \\ - 4062662189783485600 x + 145700736445997574400
\end{multline*}

\paragraph{Observations.}


As expected, with L\"{o}schian values of $k$ we found that regular hexagons are hard to beat. The smallest L\"{o}schian $k$ for which regular hexagons are beaten by another class (in this case a rectilinear one but not a semi-regular one) is $k = 112$. 

By contrast, we obtained substantial improvements on earlier results with many non-L\"{o}schian $k$. For example, we achieved $d(8)=7/5$, a considerable improvement on the value of 1.37542 reported in \cite{chy}. 


In the case of $k=156$ not only do rectilinear hexagons beat regular ones, but both are record breaking compared to the previous $d(155)$. 
We checked $f(k)$ for other cases where $k=n(n+1)$ is L\"{o}schian: $\{756, 1332,$ $2352, \dots\}$. For all subsequent $k<10^{10}$ rectilinear hexagons beat regular ones. It is natural to conjecture that this will be true for all $k>12$.


Rectilinear hexagons can usually beat semi-regular ones, especially in the case where semi-regular hexagons beat regular ones. The only $k \le 175$ for which a semi-regular hexagon beats regular ones, is not beaten by a rectilinear one, and beats all smaller $k$, are $k=8,15,33,96,99,143,168$.

Except for the few cases $k=11, 15, 18, 23, 45, 77$ and $187$, the optimal tiling has each hexagon equidistant from six same-coloured ones, i.e. $d_{01}=d_{02}=d_{12}$. Note that the exceptions completely include the class associated with the fourth degree polynomial.

For $k=80$ and $120$, a quadratic-form semi-regular $d^2(k)$ is beaten by a rectilinear example in which $d^2(k)$ is rational.


For some values of $k$ there are strikingly large values of $n$ for which $d(k)<d(k-n)$. For example, $d(29)<d(25)$, $d(77)<d(69)$, $d(98)<d(91)$, $d(187)<d(169)$.

\section{Avenues for future work}

\paragraph{More general hexagons.}

Clearly this work can be extended to the analysis of more general tilings. The stipulation that the hexagon must be rectilinear can be relaxed, for example. More ambitious would be to allow hexagons of multiple shapes, or to assign colours to sublattices that are not translations of each other. In general, we do not expect such tilings to beat the ones we describe here -- though the fact that there remain so many cases where $d(k)<d(k-1)$ may be viewed as evidence to the contrary.

\paragraph{Beyond hexagons.}

Looking beyond hexagons, however, the outlook is much rosier. We found a tiling that can be 8-coloured with $d\approx 1.444157$, a remarkable improvement over the value of $7/5$ found in this study. We also found 14- and 15-tilings that beat hexagons by similarly impressive margins: $d(14)\approx 2.260808$,  $d(15)\approx 2.346969$. Further details will be provided in a forthcoming article (see also the discussion within the Polymath16 project \cite{pol}).

\paragraph{Towards \textit{d}(6)\,=\,1.}

Finally we come to 
the only report of which we are aware that has discussed $d(k)$ in the context of tiling with identical but non-regular hexagons. Ivanov reported \cite{iva} that $d(6)$ can be as much as $\approx 0.991445$ using rectilinear hexagons. As shown in Table~\ref{sult}, we obtained a fractionally better value with the same design, $\approx 0.992076$ (see Fig.~\ref{hex1}); we do not know how Ivanov derived the parameters leading to his value. 
For comparison, the best pentagonal tiling we know has $d(6)=\sqrt{55/56}\approx 0.991031$ (
marked with an asterisk in Fig.~\ref{hex1}). 

Since a block-structured $6$-tiling of the plane excluding distance 1 can only use tiles with at most five edges, simplistic application of the methods used here cannot reach $d(6)=1$. 
However, we cannot exclude such possibility with more general polygons. The remarkable improvements to $d(8), d(14)$ and $d(15)$ noted above may inspire 
optimism in this regard.

\paragraph{Redundant colours?}

As we already noted, the nonmonotonicity of the function $d(k)$ is associated with the restrictions of the study area. 
However, if we remove all additional restrictions, 
then for some $k$ we still get $d(k)=d(k-1)$ and even $d(k)= d(k-2)$ (minimal examples with $d>0$ here are $k=5$ and $k=11$ respectively).
In other words, adding one or two colours does not change the maximum distance between tiles of the same colour. And this state of affairs is extremely surprising.

The question arises whether the function $d(k)$ is strictly monotonic or are some colours actually redundant? The available data leave this question open.

\end{document}